\documentclass{amsart}

\usepackage[latin1]{inputenc}
\usepackage[all]{xy}

\usepackage{amssymb}
\usepackage{amsthm}
\usepackage{mathtools}
\usepackage{tikz-cd}
\usepackage{color}
\usepackage{hyperref}
\usepackage{upref}
\usepackage{enumitem}
\usepackage{verbatim}

\usepackage[latin1]{inputenc}

\setlist[enumerate]{nosep}
\setlist[enumerate,1]{label={\textup{(\arabic*)}}}

\allowdisplaybreaks

\newtheorem{thm}{Theorem}[section]
\newtheorem{prop}[thm]{Proposition}
\newtheorem{cor}[thm]{Corollary}
\newtheorem{lem}[thm]{Lemma}

\theoremstyle{definition}
\newtheorem{defn}[thm]{Definition}

\newtheorem{rem}[thm]{Remark}

\numberwithin{equation}{section}

\newcommand{\secref}[1]{Section~\textup{\ref{#1}}}

\newcommand{\thmref}[1]{Theorem~\textup{\ref{#1}}}
\newcommand{\corref}[1]{Corollary~\textup{\ref{#1}}}
\newcommand{\lemref}[1]{Lemma~\textup{\ref{#1}}}
\newcommand{\propref}[1]{Proposition~\textup{\ref{#1}}}

\newcommand{\remref}[1]{Remark~\textup{\ref{#1}}}
\newcommand{\diagref}[1]{Diagram~\textup{\ref{#1}}}

\newcommand{\cc}[1]{\mathcal{#1}}

\newcommand{\MM}{\cc M}
\newcommand{\NN}{\cc N}
\renewcommand{\AA}{\cc A}
\newcommand{\BB}{\cc B}
\newcommand{\DD}{\cc D}

\newcommand{\CC}{\cc C}
\newcommand{\HH}{\cc H}

\renewcommand{\epsilon}{\varepsilon}

\newcommand{\inv}{^{-1}}
\newcommand{\wilde}{\widetilde}

\newcommand{\id}{\text{id}}

\renewcommand{\subset}{\subseteq}

\renewcommand{\bar}{\overline}

\newcommand{\variso}{\overset{\simeq}{\longrightarrow}}

\newcommand{\cst}{\ensuremath{C^*}}
\newcommand{\csta}{\ensuremath{C^*}-algebra}

\newcommand{\cstaa}{\ensuremath{C^*(\AA)}}
\newcommand{\gam}[1]{\ensuremath{\Gamma_c(\mathcal{#1})}}
\newcommand{\gac}{\Gamma_c(\AA)}
\newcommand{\gacz}{\Gamma_c(\AA^0)}
\newcommand{\gaz}[1]{\ensuremath{\Gamma_0(\mathcal{#1}^0)}}
\newcommand{\cstbb}{\ensuremath{C^*(\BB)}}

\newcommand{\cstra}{\ensuremath{C^*_r(\AA)}}
\newcommand{\cstrb}{\ensuremath{C^*_r(\BB)}}

\newcommand{\hm}{homomorphism}

\newcommand{\rep}{representation}
\newcommand{\ce}{conditional expectation}

\newcommand{\co}{\ensuremath{\mathbf{Co}}}

\newcommand{\ma}{\ensuremath{\mathbf{Max}}}

\newcommand{\midtext}[1]{\quad\text{#1}\quad}

\begin{document}

\title[Monads, comonads and equivalent subcategories]
{On groupoid-graded C*-algebras and equivalent subcategories linked via monads and comonads}
\author{Erik B\'edos}
\address{Department of Mathematics, University of Oslo, PB 1053 Blindern, 0316 Oslo, Norway}
\email{bedos@math.uio.no}

\author{S. Kaliszewski}
\address{School of Mathematical and Statistical Sciences
\\Arizona State University
\\Tempe, Arizona 85287}
\email{kaliszewski@asu.edu}

\author{John Quigg}
\address{School of Mathematical and Statistical Sciences
\\Arizona State University
\\Tempe, Arizona 85287}
\email{quigg@asu.edu}

\date{\today}


\begin{abstract}
We present a new method, involving monads and comonads from category theory,
to help establish a certain type of equivalence of subcategories.
As a case study we consider the category of topological gradings of $C^*$-algebras over 
a fixed Hausdorff \'etale groupoid. 
\end{abstract}

\maketitle

\section{Introduction}\label{intro}

Category theory has been applied to operator algebras for a long time.
Frequently this involves categories whose objects are \csta s with extra structure,
and a morphism satisfying some universal property.
(For example, crossed products of \cst-dynamical systems.)
In some cases, the universal property is tricky to prove ---
for example, with maximalization of coactions of a group $G$.
The universal property of maximalizations was proved by Fischer (see \cite{fischer}),
but after a while we noticed that there seemed to be a gap in the verification of the uniqueness aspect.
We were eventually able to supply a satisfactory proof of this uniqueness (see \cite{bkqt}),
applying what seemed to be ad hoc categorical arguments.
Quite recently we discovered that what we had in fact done was stumble across a standard construction in category theory: monads (more precisely in this instance, comonads).

In fact, the universal property in question is part of an ``adjoint situation'',
which category theorists have known for a long time is intimately connected with monads and their dual counterparts, comonads.
More precisely, the existence of a 
left
adjoint $N\dashv \inc_\NN$
to the inclusion functor $\inc_\NN$
of a full subcategory $\NN$ in a category $\CC$ involves a universal property via the natural bijections
\[
\phi_{x,y}:\CC(N x,y)\variso \CC(x,y),
\]
where
$y$ is in the subcategory defined by
$\eta_y$ being invertible.

In \secref{prelim} we explain our conventions regarding Fell bundles over \'etale groupoids,
and in \secref{monads} we do the same for monads and comonads.
Moreover, we give a useful improvement in the general theory:
removing one of the axioms in an idempotent comonad
gives an ostensibly weaker structure, which we call an \emph{idempotent semicomonad} (see the paragraph following \remref{mono}),
and then in \corref{semicomonad is comonad} we prove that the two concepts are equivalent; in \corref{semimonad is monad} we prove the dual fact for
\emph{idempotent semimonads} and \emph{idempotent monads}.

In \secref{graded} we begin the study of another example of (what we call) \emph{maximal-normal} equivalence of subcategories:
the example involves
topological gradings of \csta s over a fixed Hausdorff \'etale groupoid $G$.
A \emph{$G$-grading} is determined by a Fell bundle $\AA$ over $G$.
The ``maximal'' gradings are associated to the universal algebra \cstaa,
and the grading is \emph{topological} if it is compatible with the reduced algebra \cstra,
giving rise to the ``normal'' gradings.
In Sections~\ref{maximalization} and \ref{normalization} we introduce relevant functors $M$ and $N$,
and in \secref{max-nor} we prove that this gives another instance of the maximal-normal equivalence (see \thmref{max nor grading}).

Part of our purpose in this paper is to promulgate the view that these maximal-normal subcategories are ubiquitous in \cst-theory.
For yet another example, 
{see \cite{jeri}, where
J. Spiker (a student of the third author) is currently working on an example involving
compactifications of locally compact Hausdorff spaces and unitizations of \csta s.

\section{Preliminaries}\label{prelim}
We will use the convention that if $\phi: A \to B$ is an adjoint preserving homomorphism between some complex $*$-algebras $A$ and $B$, we will just say that $\phi$ is a homomorphism. Accordingly, any isomorphism between $A$ and $B$ is assumed to preserve adjoints, and this also applies to any representation of $A$ by bounded operators on a Hilbert space.
  
Throughout this paper, $G$ denotes an \'etale groupoid (see for instance \cite{williams2019}),
which we always assume to be Hausdorff. The unit space of $G$ is denoted by $G^0$, and $r:G\to G^0$,  $s:G\to G^0$  are the range and source maps, respectively, while $G^2:= \{ (g, h) \in G\times G : s(g)=r(h)\}$ is the set of composable pairs. When $x\in G^0$, then $G^x:=\{g\in G: r(g)=x\}$ and $G_x:=\{g\in G: s(g)=x\}$ are discrete subsets of $G$. An open subset $U$ of $G$ is called a bisection when the restrictions of $r$ and $s$ to $U$ are homeomorphisms onto $r(U)$ and $s(U)$, respectively. The family of bisections of $G$ forms then a basis for the topology on $G$.

We refer to \cite[Appendix C]{williams2007} for the definition and properties of upper semicontinuous  Banach bundles over locally compact Hausdorff spaces. By a Fell bundle $\AA$ over $G$, we will mean an upper semicontinuous Banach bundle $\AA$ over $G$ associated to a continuous, open surjection $p:\AA \to G$, equipped with a continuous multiplication map $(a, b)\mapsto ab$ from $\AA^{2} = \{ (a, b) \in \AA\times\AA: (p(a), p(b)) \in G^2\}$ into $\AA$, and a continuous involution $a\mapsto a^*$ from $\AA$ into $\AA$, satisfying the axioms spelled out by Kumjian in \cite[Definition 2.1]{kumjian}. In particular, if $x \in G^0$, then $\AA_x:= p^{-1}(x)$ is a $C^*$-algebra (with multiplication, involution and norm induced from the bundle). 
Kumjian requires the Banach bundle $\AA$ to be continuous, but upper-semicontinuity is sufficient for what follows, cf.~\cite[Section 2.2]{bussexel2012}. The restriction $\AA^0= \AA|_{G^0}$ is then an upper semicontinuous $C^*$-bundle over $G^0$.
We write $\Gamma_c(\AA)$ for the
space of continuous sections of $\AA$ with compact support, which is a $*$-algebra with respect to the convolution product and the involution given by
 \[ (f_1* f_2)(g) = \sum_{h\in G^{r(g)}} f_1(h) f_2(h^{-1}g), \quad f_1^*(g) = f_1(g^{-1})^*\]
 for $f_1, f_2 \in \Gamma_c(\AA)$ and $g\in G$.
The space $\Gamma_c(\AA^0)$ of continuous compactly supported sections of $\AA^0$ can then be canonically identified with a $*$-subalgebra of  
 $\Gamma_c(\AA)$. We note that $\Gamma_c(\AA^0)$ is a pre-$C^*$-algebra, whose  completion is the \csta\ $\Gamma_0(\AA^0)$ of sections of $\AA^0$ vanishing at infinity, equipped with the sup norm $\|\cdot\|_\infty$. We let $P_\AA: \Gamma_c(\AA)\to \Gamma_c(\AA^0)$ denote the restriction map. 
 
The space $\Gamma_c(\AA)$ is also a right inner-product $\Gamma_c(\AA^0)$-module with respect to  the right action of  $\Gamma_c(\AA^0)$ on  $\Gamma_c(\AA)$ and the  $\Gamma_c(\AA^0)$-valued inner product on  $\Gamma_c(\AA)$ given respectively by 
  \[ (f\cdot \xi)(g) = f(g) \xi(s(g))\]
 for all $f \in \Gamma_c(\AA)$, $\xi\in \Gamma_c(\AA^0)$ and  $g \in G$, 
and
    \[\langle f_1, f_2\rangle = P_\AA(f_1^**f_2),\]
    that is,
  \[  \langle f_1, f_2\rangle(x) = (f_1^** f_2)(x) = \sum_{k\in G_x} f_1(k)^*f_2(k), \]
for all $f_1, f_2 \in \Gamma_c(\AA)$ and $x \in G^0$, cf.~\cite[Proposition 3.2]{kumjian}. We denote by $L^2(\AA)$ the right Hilbert $\Gamma_0(\AA^0)$-module obtained by completing $\Gamma_c(\AA)$ with respect to the norm $\|f\|_2:= \|\langle f, f\rangle\|_\infty^{1/2}$. The left regular representation of $\Gamma_c(\AA)$, i.e., the action of $\Gamma_c(\AA)$ on itself by left multiplication, extends to an injective representation $i_{\AA, r}$ of $\Gamma_c(\AA)$ on $L^2(\AA)$ by adjointable operators.   The reduced cross-sectional $C^*$-algebra $C^*_r(\AA)$ associated to $\AA$ is by definition the closure of $i_{\AA, r}(\Gamma_c(\AA))$ with respect to operator norm.
 
We write $\cstaa$ for the  enveloping \csta\ of the $*$-algebra $\Gamma_c(\AA)$,
and $i_\AA:\Gamma_c(\AA)\to\cstaa$ for the canonical embedding.  
 For completeness, we recall 
 how it can be shown that  $C^*(\AA)$ exists (see for example \cite{bussexel2012, kranz2023}).
 As noted in \cite[p.~166]{bussexel2012}, it holds that $\|\pi(\xi)\| \leq \|\xi\|_\infty$ whenever $\pi$ is a representation of $\Gamma_c(\AA)$ on a Hilbert space and $\xi\in \Gamma_c(\AA)$ is supported in some bisection of $G$. By a partition of unity argument (for example as in the proof of \cite[Lemma 3.1.3]{Sims2017}) one deduces
  that every section in $\Gamma_c(\AA)$ is a finite sum of such $\xi$'s.  It follows then readily that  
 \[ \|f\|:= \sup\{ \|\pi(f)\|: \pi \text{ is a representation of $\Gamma_c(\AA)$ on some Hilbert space}\}\]
defines $C^*$-seminorm on $\Gamma_c(\AA)$ satisfying that $\|\xi\| \leq \|\xi\|_\infty$ whenever $\xi\in \Gamma_c(\AA)$ is supported in some bisection of $G$.  Since $i_{\AA, r}$ is injective, $\|\cdot\|$ is in fact a $C^*$-norm, and $\cstaa$ is the completion of $\Gamma_c(\AA)$ with respect to this norm.
By construction, for any \rep\ $\pi$ of \gam A on a Hilbert space $H$, there is a unique \rep\ $\wilde\pi$ of $\cst(\AA)$ on $H$ making the diagram
\[
\begin{tikzcd}
\cst(\AA) \arrow[dr,"\wilde\pi"]
\\
\gam A \arrow[u,"i_\AA"] \arrow[r,"\pi"']
&B(\HH)
\end{tikzcd}
\]
commute.
A routine use of the Gelfand--Naimark Theorem allows us to convert this to a space-free setting. For convenient reference we record this here:
\begin{lem}\label{space}
If $f\in \gam A$ then
\begin{equation}\label{new norm}
\|f\|=\sup\{\|i(f)\|:i\text{ is a \hm\ of }\gam A\text{ into a \csta}\}.
\end{equation}
Consequently, for any \hm\ $i$ of \gam A into a \csta\ $A$,
there is a unique \hm\ $\wilde i$ making the diagram
\[
\begin{tikzcd}
\cst(\AA) \arrow[dr,"\wilde i",dashed]
\\
\gam A \arrow[r,"i"'] \arrow[u,"i_\AA"]
&A
\end{tikzcd}
\]
commute, whose range is the closure of $i(\gam A)$ in $A$.
\end{lem}
\begin{proof}
Let $\|f\|'$ denote the right-hand side of \eqref{new norm}.
Then $\|\cdot\|\le \|\cdot\|'$ is clear.
For the opposite inequality, given $i$ as above, let $\pi$ be a faithful \rep\ of $A$ on a Hilbert space  $H$. Then
\[
\|i(f)\|=\|(\pi\circ i)(f)\|\le \|f\|,
\]
and taking the sup gives $\|f\|'\le \|f\|$.
\end{proof}
Applying this lemma to $i_{\AA, r}: \Gamma_c(\AA) \to C_r^*(\AA)$, we get that there is a unique homomorphism $\Lambda_\AA:= \wilde{i_{\AA, r}}: C^*(\AA)\to C_r^*(\AA)$ satisfying that  $\Lambda_\AA\circ i_\AA = i_{\AA, r}$, which is surjective.  
 
We will also need the following folklore result.

 \begin{lem} The restrictions to $\Gamma_c(\AA^0)$ of the canonical embeddings $i_\AA$ and $i_{\AA, r}$ into $\cstaa$ and $C^*_r(\AA)$, respectively, are isometric w.r.t.~$\|\cdot\|_\infty$  on $\Gamma_c(\AA^0)$. 
 \end{lem}
 \begin{proof}
  Let $\xi \in \Gamma_c(\AA^0)$. As $G^0$ is a bisection of $G$, we have
\[ \|i_{\AA,r}(\xi) \| \leq \|i_\AA(\xi)\| = \|\xi\| \leq \|\xi\|_\infty.\] 
Thus it suffices to check that $ \|\xi\|_\infty\leq\|i_{\AA,r}(\xi) \|$. Now, for every $\eta\in \Gamma_c(\AA^0)$, we have
\[\|\eta \|_2 = \|\eta^**\eta\|_\infty^{1/2} = \sup_{x\in G^0} \|\eta(x)^*\eta(x)\|^{1/2} = \|\eta\|_\infty.\]
Hence,
\begin{align*} \|i_{\AA,r}(\xi)\| &\geq \sup\{ \| \xi*\eta\|_2\mid \eta \in \Gamma_c(\AA^0), \|\eta\|_2 = 1\}\\ &= \sup\{ \| \xi \eta\|_\infty\mid \eta \in \Gamma_c(\AA^0), \|\eta\|_\infty = 1\}\\ 
&= \|\xi\|_\infty.
\end{align*}
\end{proof}

\begin{rem} \label{C-bundle} This lemma implies  that  if  $\BB$ is a $C^*$-bundle over a locally compact Hausdorff space $X$ (and we consider $X$ as an \'etale groupoid), then we have $\|f\| = \|f\|_\infty$ for every $f\in \Gamma_c(\BB)$, i.e.,  $\Gamma_0(\BB)= C^*(\BB)$. Alternatively, this could also have been deduced from \cite[Lemma 2.3]{QS}.
\end{rem}
\begin{prop} \label{Gamma_0}
Let  $\varphi:\Gamma_c(\AA^0)\to A$ be an injective \hm\ into a $C^*$-algebra $A$. Then $\varphi$ extends uniquely  to an injective \hm\ from $\Gamma_0(\AA^0)$ into $A$.
\end{prop} 
\begin{proof} 
Remark \ref{C-bundle} gives that $\|\varphi(f)\| \leq \|f\|_\infty$ for every $f \in \Gamma_c(\AA^0)$. Hence we may extend $\varphi$ to a \hm\ 
$\bar\varphi: \Gamma_0(\AA^0) \to A$. Set $K:=\ker\bar\varphi$.
It suffices to show that $K=\{0\}$.
We know that $K\cap \gam A=\{0\}$.
Take any $a\in K$. To see that we must have $a=0$,
take any $u\in G^0$, and choose $g\in C_c(G^0)$ such that $g(u)=1$.
Since $K$ and $\gam A$ are both $*$-ideals of $\Gamma_0(\AA)$, we have
$ag\in K\cap \gam A=\{0\}$, so
\[
a(u)=a(u)g(u)=(ag)(u)=0.
\]
\end{proof}

Finally, we mention that, as shown 
 in \cite[Propositions~3.6 and 3.10]{kumjian}, the restriction map $P_\AA:\Gamma_c(\AA)\to\Gamma_c(\AA^0)$
extends to a faithful conditional expectation
\[
P_\AA:\cstra\to \gaz A.
\]
Here,  $\Gamma_0(\AA^0)$ is identified with the closure of $i_{\AA, r}(\Gamma_c(\AA^0))$ in $\cstra$.

\section{Monads and comonads}\label{monads}

\subsection{}
Throughout this subsection we let $\CC$ be a category, $M:\CC\to\CC$ be a functor,  $\psi:M\to \id_{\CC}$ be a natural transformation, and  $\CC^\psi$ be the full subcategory of $\CC$ whose objects $y$ satisfy that $\psi_y:My \to y\text{ is an isomorphism}$. When $\DD$ is a full subcategory of $\CC$, then $ {\rm Inc}_{\DD}:\DD\to \CC$ denotes the inclusion functor.

To prove the universal property of maximalization of coactions, we used in 
\cite[Proposition~4.1]{bkqt}
the following result:

\begin{prop}\label{aup}
Assume that $M\psi:M^2\to M$ is a natural isomorphism.
If 
$y$ is an object in $\CC^\psi$ and $f:y\to x$  is a morphism in $\CC$, then there is a unique morphism 
$\wilde f$ $:y\to Mx$ 
in $\CC$ making the diagram
\[
\begin{tikzcd}
y \arrow[r,"\wilde f",dashed] \arrow[dr,"f"']
&Mx \arrow[d,"\psi_x"]
\\
&x
\end{tikzcd}
\]
commute.
\end{prop}

Existence is easy: $\wilde f:=Mf\circ \psi_y\inv$.
For uniqueness, if $h,k$ both factorize $f$,
apply $M$ to the factorizations
and invertibility of $M\psi_x$
to get $Mh=Mk$.
Then naturality of $\psi$ applied to $h,k:y\to Mx$
gives
$h=\psi_{Mx}\circ Mh\circ \psi_y\inv
=\psi_{Mx}\circ Mk\circ \psi_y\inv
=k$.

Note that the condition that $M\psi$ is a natural isomorphism really only rests upon whether $(M\psi)_x$ is an isomorphism for every object $x$ of $\CC$;
$M\psi$ is of course natural,
since $M$ is a functor and $\psi$ is a natural transformation.

 \begin{rem} One can also consider the condition that $\psi M:M^2\to M$ is a natural isomorphism\footnote{As above, note that the only issue here is whether $(\psi M)_x$ is an isomorphism for every object $x$ of $\CC$.}, which says that $Mx$ is an object of $\CC^\psi$ for every object $x$ of $\CC$. When this condition holds, we will let $M^\psi: \CC \to \CC^\psi$ denote the functor obtained by  considering $M$ as a functor from $\CC$ to $\CC^\psi$; we then have $M = {\rm Inc}_{\CC^\psi}\circ M^\psi$.    We note that in general we have $M\psi \not= \psi M$, so this condition should not be confused with the one assumed in Proposition \ref{aup}. 
 \end{rem}
 
 Next, we recall that the pair  $(M, \psi)$ is  called an \emph{idempotent comonad} on $\CC$ when 
 $M\psi = \psi M$ and $M\psi$ is a natural isomorphism.\footnote{This is not the most common definition, but the one that fits best with our discussion. The connection with the original definition, which is related to the more general concept of comonad in category theory, is explained in Grandi's book (\cite[page 141]{grandis}).}  
  
    \begin{rem} \label{mono} Assume that  $\psi_x$ is a monomorphism for every object $x$ of $\CC$. 
  Then $\psi M = M \psi$. Thus the following conditions are equivalent in this case:
 
 \begin{itemize} 
 \item[$(a)$] $M\psi:M^2\to M$ is a natural isomorphism. 
 \item[$(b)$] $\psi M:M^2\to M$ is a natural isomorphism. 
  \item[$(c)$] $(M, \psi)$ is an idempotent comonad on $\CC$.
 \end{itemize}

Indeed, let $x$ be an object of $\CC$. Then, by naturality of $\psi$, we get that \[\psi_x \circ \psi_{Mx}= \psi_x \circ (M\psi_x).\] 
Since $\psi_x$ is a monomorphism, this implies that $\psi_{Mx}= M\psi_x$. Thus, $\psi M = M \psi$, as asserted.
\end{rem}

Let us say that  $(M, \psi)$ is an \emph{idempotent semicomonad} on $\CC$ when $M\psi$ and  $\psi M$ are both natural isomorphisms (without requiring that $M\psi = \psi M$). Obviously, an idempotent comonad is an idempotent semicomonad. 
We will soon see that the converse holds.

We pause to interpret the above in the context of
our paper 
\cite{bkqt}:
we had
a category \co\ of coactions of a fixed locally compact group,
a functor \ma, and a natural transformation $\psi:\ma\to \id_{\co}$.
It follows readily from 
\cite{bkqt}
that $(\ma, \psi)$ gives an example of an idempotent semicomonad on \co.
 Our interest in the concept of idempotent semicomonad comes from the following result:

 \begin{prop} \label{corefl}  Assume that  $(M, \psi)$ is an idempotent semicomonad on $\CC$.

  Then ${\rm Inc}_{\CC^\psi}: \CC^\psi \to \CC$ is right-adjointable and $M^\psi$ is a right adjoint. Thus  $\CC^\psi$ is a coreflective subcategory of $\CC$. Moreover, $\psi$ is the counit of the adjunction ${\rm Inc}_{\CC^\psi} \dashv M^\psi$.
 \end{prop} 
 
 \begin{proof} 
 Proposition \ref{aup} gives that the pair $ (M^\psi x, \psi_x) = (Mx, \psi_x)$ is a universal morphism from $\CC^\psi$ to $x$ for every object $x$ of $\CC$. Hence the assertions follow from general principles (see for example 
\cite[Section~2]{reflective}).
 \end{proof} 
 On the other hand, it is well-known that an adjunction of functors gives rise to a comonad (resp.~monad), and the fact that this comonad (resp.~monad) is idempotent can be characterized  in different ways, see for example  
\cite[Proposition~3.8.2]{grandis}.
 Applying this result to the comonad of the adjunction associated to a coreflective subcategory, we get
  \begin{prop}\label{corefl-idem}
  Assume that $\MM$ is a coreflective subcategory of $\CC$ and let  $ M^\sharp: \CC\to \MM$ denote a right-adjoint of  ${\rm Inc}_{\MM}$. Set $M:=  {\rm Inc}_{\MM}\circ M^\sharp : \CC \to \CC$ 
 and let $\psi : M \to {\rm id}_\CC$ denote the counit of the adjunction ${\rm Inc}_{\MM} \dashv M^\sharp$.  
 Then the following conditions are equivalent:
 \begin{itemize}
   \item $(M, \psi)$ is an idempotent comonad on $\CC$.
   \item $\psi_y: M y \to y$ is an isomorphism 
   for every object $y$ of $\MM$. 
    \item $M\psi : M^2 \to M$ is a natural isomorphism. 
  \item $\psi M : M^2 \to M$ is a natural isomorphism. 
  \end{itemize}
  \end{prop}
\begin{proof} This follows from the equivalence between (i), (iii$^*$), (iv$^*$), and (v$^*$) in 
\cite[Proposition 3.8.2]{grandis}.
\end{proof} 
 Combining Propositions \ref{corefl} and \ref{corefl-idem} gives:
 
 \begin{cor}\label{semicomonad is comonad}
 Assume that  $ (M, \psi)$ is an idempotent semicomonad on $\CC$.
 Then  $(M, \psi)$ is an idempotent comonad on $\CC$.
   \end{cor}
 In particular, we get  that $(\ma, \psi)$ is an idempotent comonad on $\co$, i.e., $\ma\,\psi =\psi \,\ma$ (and $\ma\,\psi$ is an isomorphism).
 
\begin{rem}\label{mono-max} In Remark \ref{mono}, we saw that if
$\psi_x$ is a monomorphism for every object $x$ of $\CC$, then $M\psi=\psi M$.
We don't know whether it has been observed before that the monomorphicity property holds in the case of  $(\ma, \psi)$.
To see this, let's operate abstractly. Assume that $(M,\psi)$ satisfies the universal property of Proposition \ref{aup} 
and that the components $\psi_x$ are epimorphisms for every object $x$ of $\CC$. (Note that these conditions are satisfied for $(\ma, \psi)$).

Let $x$ be an object of $\CC$. We want to show that $\psi_x$ is a monomorphism. Consider $\phi,\tau:y\to Mx$ satisfying $\psi_x\phi=\psi_x\tau$.
First, if $y\in \CC^\psi$, then $\phi=\tau$ by the uniqueness in the universal property.
Then in general we have
$\psi_x\phi\psi_y=\psi_x\tau\psi_y$, so $\phi\psi_y=\tau\psi_y$ by the preceding.
But then $\phi=\tau$ since $\psi_y$ is an epimorphism. 

This gives a different way to deduce that that $\ma\, \psi=\psi\,\ma$. 
To our knowledge this particular fact has not been noted in the literature.
\end{rem}

\begin{rem} \label{coalgebraic} Let $(M, \psi)$ be an idempotent comonad on $\CC$. Then $\CC^\psi$ coincides with the full subcategory $\mathbf{Coalg}(M)$ of $\CC$ whose objects are the so-called \emph{coalgebraic} objects of $\CC$,
that may be characterized in several ways, cf.~\cite[3.8.5]{grandis}.  
For instance, it follows that $x$ is an object of $\CC^\psi$ if and only if $x$ is isomorphic to  $Mx$ in $\CC$.
\end{rem} 
 
 \subsection{}\label{monad subsec}
 For completeness we also consider the dual situation. Throughout this subsection, we let $N:\CC\to\CC$ be a functor, $\eta:\id_\CC\to N$ be a natural transformation, and  $\CC^\eta$ be the full subcategory of $\CC$ whose objects $y$  satisfy that $\eta_y\text{ is an isomorphism}$.
 
 The dual version of Proposition \ref{aup}, which can be used to show the universal property of normalization of coactions, is as follows.
\begin{prop}\label{aup2} 
Asssume that  $N\eta:N\to N^2$
is a natural isomorphism. 
 
If  $y$ is an object of $\CC^\eta$ and
$g:x\to y$ is a morphism  in $\CC$, then there is a unique morphism 
$g': Nx\to y$ in $\CC$  such that $g' \circ \eta_x = g'$.
\end{prop}

We also state the dual versions of  the results in the previous subsection. 
We recall that 
 the pair $(N, \eta)$ is called an \emph{idempotent monad} on $\CC$ when 
 $N\eta = \eta N$ and $N\eta$ is a natural isomorphism. 
   \begin{rem} \label{epi} Assume that  $\eta_x$ is an epimorphism for every object $x$ of $\CC$. 
  Then $\eta N = N \eta$. Thus the following conditions are equivalent in this case:
 
 \begin{itemize} 
 \item[$(a)$] $N\eta:N\to N^2$ is a natural isomorphism. 
 \item[$(b)$] $\eta N:N\to N^2$ is a natural isomorphism. 
  \item[$(c)$] $(N, \eta)$ is an idempotent monad on $\CC$.
 \end{itemize}
 \end{rem}
 
 When $\eta N:N\to N^2$ is a natural isomorphism, that is,  $Nx$ is an object of $\CC^\eta$ for every object $x$ of $\CC$,  we will let $N^\eta: \CC \to \CC^\eta$ denote the functor obtained by  considering $N$ as a functor from $\CC$ to $\CC^\eta$; thus $N = {\rm Inc}_{\CC^\eta}\circ N^\eta$. 
 
 We will say that  $(N, \eta)$ is an \emph{idempotent semimonad} on $\CC$ when $N\eta$ and  $\eta N$ are both natural isomorphisms (without requiring that $N\eta = \eta N$). 
 An idempotent monad is clearly an idempotent semimonad.  
 
 It is easy to see from  
\cite{bkqt}
 that $(\mathbf{Nor}, \eta)$ is an idempotent semimonad on $\co$. 

  \begin{prop} \label{refl} Assume that $(N, \eta)$ is an idempotent semimonad on $\CC$.
  Then ${\rm Inc}_{\CC^\eta}: \CC^\eta \to \CC$ is left-adjointable and $N^\eta$ is a left adjoint. Thus  $\CC^\eta$ is a reflective subcategory of $\CC$. Moreover, $\eta$ is the unit of the adjunction $N^\eta \dashv {\rm Inc}_{\CC^\eta}$.
\end{prop} 

\begin{prop} Assume that $\NN$ is a reflective subcategory of $\CC$ and let  $ N^\sharp: \CC\to \NN$ denote a left-adjoint of  ${\rm Inc}_{\NN}$. Set $N:=  {\rm Inc}_{\NN}\circ N^\sharp : \CC \to \CC$ 
 and let $\eta : 1_\CC\to N$ be the unit of the adjunction $N^\sharp \dashv {\rm Inc}_{\NN}$.  
Then the following conditions are equivalent: 
  \begin{itemize}
   \item  $ (N, \eta)$ is an idempotent monad. 
     \item $\eta_y: y \to N y$ is an isomorphism 
      for every object $y$ of $\NN$. 
    \item $N\eta : N \to N^2$ is a natural isomorphism. 
  \item $\eta N : N \to N^2$ is a natural isomorphism. 

\end{itemize}
\end{prop}
 \begin{cor}\label{semimonad is monad}
 Assume that  $ (N, \eta)$ is an idempotent semimonad on $\CC$.
 Then 
   $(N, \eta)$ is an idempotent monad on $\CC$.
   \end{cor}
   
 In particular, we get  that $(\mathbf{Nor}, \eta)$ is an idempotent monad on $\co$, i.e., $\mathbf{Nor} \,\eta =\eta \,\mathbf{Nor}$ (and $\mathbf{Nor}\,\eta$ is an isomorphism).
Alternatively, arguing  analogously as in Remark \ref{mono-max}, one can show that $\eta_x$ is an epimorphism for every object $x$ in $\co$, and use Remark \ref{epi} to deduce that $\mathbf{Nor} \,\eta =\eta \,\mathbf{Nor}$.

\begin{rem} \label{algebraic} Let $(N, \eta)$ be an idempotent monad on $\CC$. Then $\CC^\eta$ coincides with the full subcategory 
$\mathbf{Alg}(N)$ 
of $\CC$ whose objects are the so-called \emph{algebraic} objects of $\CC$,
that 
may be characterized in several ways, 
cf.~\cite[3.8.3]{grandis}.  
For instance, it follows that $x$ is an object of $\CC^\eta$ if and only if $x$ is isomorphic to  $Nx$ in $\CC$.
\end{rem} 
 
\subsection{} \label{MandN}
Let us
connect our considerations in the two preceding subsections with 
our works
\cite{reflective, asspairs} on reflective-coreflective equivalence.  
The setup 
will
now be 
as follows.  
 
\begin{itemize}
\item $\CC$ is a category,
 \item $M:\CC\to\CC$ and $N:\CC\to\CC$ are functors, 
 
 \item  $\psi:M\to \id_{\CC}$ and $\eta:\id_\CC\to N$ are natural transformations such that 
 $M\psi, \psi M : M^2 \to M$ and  $N\eta, \eta N: N\to N^2$ are all natural isomorphisms.
 \end{itemize}
 Note that this is equivalent to requiring  that  $(M, \psi)$ is an idempotent semicomonad on $\CC$ and $(N, \eta)$ is an idempotent semimonad on $\CC$.  
By Corollary \ref{semicomonad is comonad} and Corollary \ref{semimonad is monad}, this implies that  $(M, \psi)$ is an idempotent comonad on $\CC$ and $(N, \eta)$ is an idempotent monad on $\CC$.

We set $\MM:=\CC^\psi$ and $ \NN:= \CC^\eta$. As in the two previous subsections we let  $M^\psi$ denote $M$ considered as a functor from $\CC$ to $\MM$,  and $N^\eta$ denote $N$ considered as a functor from $\CC$ to $\NN$. Note that for every morphism $f$ in $\CC$, we have $M^\psi f= Mf$ and $N^\eta f = Nf$. The picture is then
\[
\xymatrix@C+30pt{
\MM \ar@<1ex>[r]^{\inc_\MM} \ar@/^2pc/[rr]^{N^\eta \circ \,{\rm Inc}_{\MM}}
&\CC \ar@<1ex>[r]^{N^\eta} \ar@<1ex>[l]^{M^\psi}
&\NN \ar@<1ex>[l]^{\inc_\NN} \ar@/^2pc/[ll]^{M^\psi \circ \,{\rm Inc}_{\NN}}
}
\]
and we have that
 \begin{itemize}
 \item $\MM$ is a full coreflective subcategory of $\CC$,
 \item $\psi$ is the counit of the adjunction ${\rm Inc}_{\MM} \dashv M^\psi$,
 \item $ \NN$ is a full reflective subcategory of $\CC$, 
 \item  $\eta$ is the unit of the adjunction $N^\eta \dashv {\rm Inc}_{\NN}$,
 \item $N^\eta \circ {\rm Inc}_{\MM} \dashv M^\psi \circ {\rm Inc}_{\NN}$.
 \end{itemize}

Our reflective-coreflective type equivalence between $\MM$ and $\NN$, 
\cite[Corollary 4.4]{reflective}, which says that the adjunction $N^\eta \circ {\rm Inc}_{\MM} \dashv M^\psi \circ {\rm Inc}_{\NN}$ is an adjoint equivalence  between $\MM$ and $\NN$,  
will then be
 satisfied whenever the following two conditions hold:
\begin{itemize}
 \item[a)] for each object $x$ in $\NN$, the pair $(x,\psi_x)$ is an
initial object in the comma category $Mx\downarrow\NN$.
\item[b)] 
for
each object $x$ in $\MM$, the pair $(x,\eta_x)$ is a final object in the comma
category $\MM\downarrow Nx$. 
 \end{itemize}
(See 
\cite[Theorem 4.3]{reflective} for some equivalent formulations of these conditions).

Moreover, if we, instead of a) and b), make the stronger assumption:
\begin{itemize} 
\item[c)] 
for
each object $x$ in $\CC$, $(Nx, \eta_x\circ \psi_x)$ is an initial object in $Mx\downarrow\NN$ and $(Mx, \eta_x\circ \psi_x)$ is a final object in $\MM\downarrow Nx$, 
 \end{itemize}
 or, equivalently,
 \begin{itemize} 
\item[c')]   the natural transformations $M\eta: M\to MN$ and $N\psi: NM \to N$ are both isomorphisms,
 \end{itemize}
then 
we get that 
the adjunction $N^\eta \circ {\rm Inc}_{\MM} \dashv M^\psi \circ {\rm Inc}_{\NN}$ is a ``maximal-normal" adjoint equivalence  between $\MM$ and $\NN$.
Equivalently, cf.~\cite[Theorem 2.2]{asspairs}, we get that $\MM$ and $\NN$ form an \emph{associated pair} in the sense of \cite{KL}.  As explained in \cite{asspairs}, this can most easily be described by saying that for every morphism $f$ in $\CC$, $Mf$ is an isomorphism if and only $Nf$ is an isomorphism.  

The fact that condition c') is equivalent to 
condition c)
 follows from \cite[Theorem 3.4]{reflective} (see the proof of \cite[Proposition 2.1]{asspairs}). We also mention that \cite[Proposition 5.3 and Corollary 5.4]{reflective} then gives that 
\[N^\eta \simeq (N^\eta \circ {\rm Inc}_{\MM}) \circ M^\psi \text{ and  } M^\psi \simeq (M^\psi \circ {\rm Inc}_{\NN}) \circ N^\eta.\]

\section{Graded C*-algebras}\label{category}\label{graded}

As a case study, 
it is natural to
consider 
categories whose objects are $C^*$-algebras that are graded over a fixed discrete group $G$. 
We recall from Exel's book \cite{exelbook} that a $C^*$-algebra $A$ is said to be 
\emph{$G$-graded} when there exists a linearly independent family $\{A_g\}_{g\in G}$ of norm-closed subspaces of $A$, called the \emph{grading family of $A$},  satisfying that  the algebraic direct sum 
 of the $A_g$'s is norm-dense in $A$ and 
\[ A_g A_h \subseteq A_{gh}, \quad  (A_{g})^{*} = A_{g^{-1}} \]
for all $g,h \in G$.

We will actually generalize somewhat from groups: we allow $G$ to be a Hausdorff  \'etale groupoid.
First of all,
if $\AA$ and $\BB$ are Fell bundles over $G$,
we define
a \emph{morphism} $\varphi:\AA\to \BB$ to be
a continuous map $\varphi:\AA\to \BB$
that restricts to linear maps $\varphi_x:A_x\to B_x$ for all $x\in G$
and is multiplicative and involutive in the sense that
\begin{align*}
\varphi(ab)=\varphi(a)\varphi(b)
\midtext{and}
\varphi(a^*)=\varphi(a)^*.
\end{align*}
(Of course, in the first we require $a,b$ to be in fibers over composable elements of $G$).
It is then folklore (and in any case is a routine exercise to check) that the following result holds.

\begin{lem}\label{integrate}
For every Fell-bundle morphism $\varphi:\AA\to\BB$ there is a unique
\hm\ $\varphi^c:\gam A\to \gam B$
given by
\[
\varphi^c(f)=\varphi\circ f.
\]
Moreover, the assignments
$\AA\mapsto \gam A$ on objects
and
$\varphi\mapsto \varphi^c$ on morphisms define a functor from the category of Fell bundles over $G$ with morphisms to the category of $*$-algebras with \hm s. 
\end{lem}

To every Fell-bundle morphism $\varphi:\AA\to \BB$ one can also associate a canonical \hm\ $\varphi^0:\Gamma_0(\AA^0) \to \Gamma_0(\BB^0)$. Indeed, it is clear that the restriction $\varphi^0$ of $\varphi^c$ to $\Gamma_c(\AA^0)$  is a \hm\ into $\Gamma_c(\BB^0)$. Considering  $\varphi^0$ as taking values in $\Gamma_0(\BB^0)$ (see \remref{C-bundle}), 
we can simply apply \propref{Gamma_0}.

\begin{prop}
For every Fell-bundle morphism $\varphi:\AA\to \BB$
there is a unique \hm\ $\varphi^m:C^*(\AA)\to C^*(\BB)$
making the diagram
\[
\begin{tikzcd}
\cstaa \arrow[r,"\varphi^m",dashed]
&\cstbb
\\
\gam A \arrow[r,"\varphi^c"'] \arrow[u,"i_\AA"]
&\gam B \arrow[u,"i_\BB"']
\end{tikzcd}
\]
commute.
\end{prop}

\begin{proof}
Apply \lemref{space} with
$i=i_\BB\circ\varphi^c$.
\end{proof}

By a  \emph{$G$-grading} of a $C^*$-algebra $A$ we will mean a pair $(\AA, i)$ where $\AA$ is a Fell bundle over $G$ and 
$i:\Gamma_c(\AA)\to A$ is an injective \hm\  with dense range.
We note that \lemref{space} tells us that
there is a unique \hm\ $\psi_A:C^*(\AA) \to A$
such that
\[
\psi_A\circ i_\AA=i,
\]
which is surjective. We will 
later frequently use the more precise notation $\psi_{(A, \AA, i)}$ instead of $\psi_A$. 

We also note that Corollary \ref{Gamma_0} gives that the restriction $i_0$ of $i$ to $\Gamma_c(\AA^0)$ extends uniquely to a isomorphism from $\Gamma_0(\AA^0)$ onto 
the closure of $i(\Gamma_c(\AA^0))$ in $A$. In the sequel, when $(\AA, i)$ is a $G$-grading of $A$, we will consider $\Gamma_0(\AA^0)$ as a $C^*$-subalgebra of $A$ using this identification. 

A $G$-grading $(\AA, i)$ of \csta\ $A$ is called \emph{topological} if
there is a (necessarily unique) \ce\ $E_A:A\to \Gamma_0(\AA^0)$
making the diagram
\[
\begin{tikzcd}
\Gamma_c(\AA) \arrow[r,"i"] \arrow[d,"P_\AA"']
&A \arrow[dl,"E_A",dashed]
\\
\Gamma_0(\AA^0)
\end{tikzcd}
\]
commute. 
Thus, restricting to the image of $\gam A$, we are asking that the map
$f\mapsto f|_{G^0}$ on \gam A (more precisely, on the image of \gam A in $A$ under the map $i:\gam A\to A$), is bounded in the $A$-norm.
We will later frequently 
use the more precise notation $E_{(A, \AA, i)}$ instead of $E_A$.

Throughout, for clarity, we restrict our attention to topological gradings, although some of our results do not need that extra assumption.

\begin{rem}
 Let $\AA$ be a Fell bundle over $G$. Then  $(\AA, i_{\AA, r})$ is a topological $G$-grading of $C_r^*(\AA)$ (with $E_{C_r^*(\AA)} = P_\AA$), and  $(\AA, i_{\AA})$ is a topological $G$-grading of $C^*(\AA)$ (with $E_{C^*(\AA)} = P_\AA\circ \Lambda_\AA$).
\end{rem} 

\begin{prop}\label{eta prop}
For every topological $G$-grading $(\AA, i)$ of a $C^*$-algebra $A$
there is a unique \hm\ $\eta_A$ making the diagram
\[
\begin{tikzcd}
\cstaa \arrow[dr,"\psi_A"] \arrow[dd,"\Lambda_\AA"']
\\
&A \arrow[dl,"\eta_A",dashed]
\\
\cstra
\end{tikzcd}
\]
commute.
Moreover, $\eta_A$ is also the unique \hm\  making the diagram 
\[
\begin{tikzcd}
\Gamma_c(\AA) \arrow[r,"i"] \arrow[d,"i_{\AA, r}"']
&A \arrow[dl,"\eta_A"]
\\
\cstra  
\end{tikzcd}
\]
commute. We will later frequently
use the more precise notation $\eta_{(A, \AA, i)}$ instead of $\eta_A$.
\end{prop}

\begin{proof}
Enlarge the first diagram:
\begin{equation}\label{eta}
\begin{tikzcd}
\cstaa \arrow[dr,"\psi_A"] \arrow[dd,"\Lambda_\AA"']
\\
&A \arrow[dl,"\eta_A",dashed] \arrow[dd,"E_A"]
\\
\cstra \arrow[dr,"P_\AA"']
\\
&\gaz A,
\end{tikzcd}
\end{equation}
where $E_A$ is the conditional expectation in the definition of topological grading.
The outer parallelepiped commutes when at the top we restrict to
$i_\AA(\Gamma_c(\AA))$,
and hence commutes at the \cst-level by density and continuity.

To show the existence of a  \hm\ $\eta_A: A \to \cstra$ such that $\eta_A\circ \psi_A = \Lambda_\AA$, it suffices to show that
\[
\ker \psi_A\subset \ker \Lambda_\AA.
\]
The important point is that the conditional expectation 
$P_\AA$ is faithful.
Thus, if $a\in \ker \psi_A$, then so is $a^*a$, hence
\[ P_\AA(\Lambda_\AA(a^*a)) = E_A(\psi_A(a^*a) = 0,\]
and we conclude that $\Lambda_\AA(a^*a)=0$, and so $a\in \ker \Lambda_\AA$, as desired. 
 
 It then follows that
 \[ \eta_A\circ i = \eta_A \circ \psi_A \circ i_\AA = \Lambda_\AA\circ i = i_{\AA,r},\]
 showing that the second diagram commutes. If $\eta: A \to \cstra$ is  a \hm\ satisfying that $\eta\circ i =  i_{\AA,r}$, then $\eta$ and $\eta_A$ agree on $i(\Gamma_c(\AA))$, hence on $A$ by density and continuity. Finally, if $\eta: A \to \cstra$ is a \hm\ satisfying that $\eta\circ \psi_A = \Lambda_\AA$, then we get as above that $\eta\circ i =  i_{\AA,r}$, which implies that $\eta=\eta_A$. 
\end{proof}

\begin{lem}\label{catone}
There is a category 
$\CC^{\rm top}_G$
in which the objects are 
triples $(A,\AA,i)$ where 
$(\AA, i) $ is a topological $G$-grading of a \csta\ $A$,
and in which a morphism $(\phi,\varphi):(A,\AA,i)\to (B,\BB,j)$ is a
\hm\ $\phi:A\to B$ together with a
Fell-bundle morphism $\varphi:\AA\to \BB$
such that
the diagram
\[
\begin{tikzcd}
\gam A \arrow[r,"\varphi^c"] \arrow[d,"i"']
&\gam B \arrow[d,"j"]
\\
A \arrow[r,"\phi"']
&B
\end{tikzcd}
\]
commutes.
\end{lem}

\begin{proof}
Most of this is obvious,
and a moment's consideration of the diagram
\[
\begin{tikzcd}
\gam A \arrow[r,"\varphi^c"] \arrow[dd,"i"']
&\gam B \arrow[dr,"\gamma^c"] \arrow[dd,"j" near end]
\\
&&\gam D \arrow[dd,"k"]
\arrow[from=llu,"\gamma^c\circ\varphi^c"' near start,crossing over]
\\
A \arrow[r,"\phi"] \arrow[drr,"\varsigma\circ\phi"']
&B \arrow[dr,"\varsigma"]
\\
&&D
\end{tikzcd}
\]
shows how the composition $(\varsigma,\gamma)\circ (\phi,\varphi)$ of morphisms works.
\end{proof}
It is easy to see that the morphisms in $\CC^{\rm top}_G$  respect the conditional expectations:
\begin{lem}\label{morph-cond-exp}
For any morphism $(\phi,\varphi):(A,\AA,i)\to (B,\BB,j)$ in $\CC^{\rm top}_G$ we have a commutative diagram
\begin{equation}\label{diag}
\begin{tikzcd}
A \arrow[r,"\phi"] \arrow[d,"E_{(A,\AA,i)}"']
&B \arrow[d,"E_{(B,\BB,j)}"]
\\
\gaz A \arrow[r,"\varphi^0"']
&\gaz B
\end{tikzcd}
\end{equation}
where 
$\varphi^0:\gaz A \to \gaz B$ 
is the homomorphism canonically associated to $\varphi$.
\end{lem}
\begin{proof}
\diagref{diag} sits inside a larger one:
\[
\begin{tikzcd}
\gac \arrow[rrr,"\varphi^c"] \arrow[dr,"i"] \arrow[dd,"P_\AA"']
&&&\gam B \arrow[dl,"j"'] \arrow[dd,"P_\BB"]
\\
&A \arrow[r,"\phi"] \arrow[dl,"E_A"']
&B \arrow[dr,"E_B"]
\\
\gacz \arrow[rrr,"\varphi^0"']
&&&\gaz B,
\end{tikzcd}
\]
where we simplified the notation a bit for this proof.
The upper quadrilateral,
the two triangles,
and the outer rectangle
all commute, and then a routine diagram chase,
using density of the image if $i$,
shows that
the lower quadrilateral commutes, as desired.
\end{proof}

\section{The maximalization functor}\label{maximalization}

In this section we will apply our comonad technique,
and for this we will need:
\begin{enumerate}
\item
a functor $M:\CC^{\rm top}_G\to \CC^{\rm top}_G$,
and

\item
a natural transformation $\psi:M\to \id_{\CC^{\rm top}_G}$
\end{enumerate}
such that
$M\psi$ and $\psi M$ are both natural isomorphisms from $M^2$ to $M$.

\begin{prop}\label{max}
There is a functor $M:\CC^{\rm top}_G\to \CC^{\rm top}_G$ given on objects by
\[
M(A,\AA,i)=(\cstaa,\AA,i_\AA)
\]
and on morphisms $(\phi,\varphi):(A,\AA,i)\to (B,\BB,j)$ by
\[
M(\phi,\varphi)=(\varphi^m,\varphi):(\cstaa,\AA,i_\AA)\to (\cstbb,\BB,i_\BB).
\]
\end{prop}
\begin{proof} Straightforward.
\end{proof}

We now promote our \hm s $\psi_{(A,\AA,i)}:\cstaa\to A$ 
to a natural transformation $\psi:M\to \id_{\CC^{\rm top}_G}$.
Actually, we must have a morphism in $\CC^{\rm top}_G$,
so we embellish
\[
(\psi_{(A,\AA,i)},\id_\AA):M(A,\AA,i)\to (A,\AA,i)
\]
and identify $\psi_{(A,\AA,i)}$ with $(\psi_{(A,\AA,i)},\id_\AA)$.

\begin{prop}\label{psi natural}
The \hm s $\psi_{(A,\AA,i)}$ combine to give a natural transformation
\[
\psi:M\to \id_{\CC^{\rm top}_G}.
\]
\end{prop}

\begin{proof}
This follows immediately from commutativity of the diagram
\[
\begin{tikzcd}
(\cstaa,\AA,i_\AA) 
\arrow[r,"{(\varphi^m,\varphi)}"] 
\arrow[d,"\psi_{(A,\AA,i)}"']
&(\cstbb,\BB,i_\BB) \arrow[d,"\psi_{(B,\BB,j)}"]
\\
(A,\AA,i) \arrow[r,"{(\phi,\varphi)}"']
&(B,\BB,j).
\end{tikzcd}
\]
\end{proof}

\begin{defn}
We denote by $\MM_G$
the full subcategory of $\CC^{\rm top}_G$ whose objects $(A,\AA,i)$
satisfy the property that
$\psi_{(\cstaa,\AA,i_\AA)}:(\cstaa,\AA,i_\AA)\to (A,\AA,i)$ is an isomorphism. In other words,  $\MM_G = \big(\CC^{\rm top}_G\big)^\psi$. 
\end{defn}

Now we investigate the compositions $M\psi$ and $\psi M$.

\begin{prop}\label{M psi}
$M\psi:M^2\to M$ is a natural isomorphism
\end{prop}

\begin{proof}
We begin by recording the formulas for
$M^2$ and $M\psi$.

First, we claim that in fact $M^2=M$.
Indeed, for an object $(A,\AA,i)$,
we have
\[
M^2(A,\AA,i)
=M(\cstaa,\AA,i_\AA)
=(\cstaa,\AA,i_\AA)
=M(A,\AA,i).
\]
For a morphism $(\phi,\varphi):(A,\AA,i)\to (B,\BB,j)$
we have
\[
M(\phi,\varphi)=(\varphi^m,\varphi):(\cstaa,\AA,i_\AA)\to (\cstbb,\BB,i_\BB)
\]
and maximalizing this does nothing, so
\begin{align*}
M^2(\phi,\varphi)
&=M(\varphi^m,\varphi)
=(\varphi^m,\varphi)
=M(\phi,\varphi),
\end{align*}
proving the claim.

We turn to $M\psi$.
For an object $(A,\AA,i)$,
we have
\[\psi_{(A,\AA,i)}=(\psi_{(A,\AA,i)},\id_\AA):(\cstaa,\AA,i_\AA)\to (A,\AA,i),\]
so
\begin{align*}
M\psi_{(A,\AA,i)}
&=((\id_\AA)^m,\id_\AA)
=(\id_{C^*(\AA)},\id_\AA)
\\&=\id_{(\cstaa,\AA,i_\AA)}
=\id_{M(A,\AA,i)}.
\end{align*}
Therefore we get
$M\psi=\id:M\to M$, which is certainly a natural isomorphism.
\end{proof}

\begin{rem}
In \propref{M psi} above we saw that $M\psi$ is in fact the identity natural isomorphism on the functor $M$.
The underlying reason why this turned out to be so trivial
is because we kept the Fell bundle $\AA$ as part of the objects in $\CC^{\rm top}_G$,
and this is what lead to $M^2=M$,
and then it was not surprising that the resulting natural isomorphism on $M$ turned out to be the identity.
\end{rem}

\begin{prop}\label{psi M}
$\psi M:M^2\to M$ is a natural isomorphism.
\end{prop}

\begin{proof}
This is similar to, and easier than, \propref{M psi}.
For an object $(A,\AA,i)$,
we have $M(A,\AA,i)=(\cstaa,\AA,i_\AA)$,
and
in \lemref{space} the relevant diagram is
\[
\begin{tikzcd}
\cstaa \arrow[dr,"\psi_{\cstaa}"]
\\
\gam A \arrow[u,"i_\AA"] \arrow[r,"i_\AA"']
&\cstaa,
\end{tikzcd}
\]
so by uniqueness we must have $\psi_{\cstaa}=\id_{C^*(\AA)}$.
Therefore we have
\[
(\psi M)_{(A,\AA,i)}=\psi_{M(A,\AA,i)}=(\id_{C^*(\AA)},\id_{\AA})=\id_{(C^*(\AA),\AA,i_\AA)},
\]
which is an isomorphism. Thus  $\psi M$ is a natural isomorphism.
\end{proof}

We now immediately conclude the following from \propref{corefl}:

\begin{thm}
Let $M:\CC^{\rm top}_G\to \CC^{\rm top}_G$ and $\psi:M\to \id_{\CC^{\rm top}_G}$ be as above.
Then:
\begin{enumerate}
\item
$(M,\psi)$ is an idempotent comonad on $\CC^{\rm top}_G$.

\item
$\inc_{\MM_G}:\MM_G\to \CC^{\rm top}_G$ is right-adjointable and
$M^\psi$ is a right adjoint.

\item
$\MM_G$ is a coreflective subcategory of $\CC^{\rm top}_G$.

\item
$\psi$ is the counit of the adjunction $\inc_{\MM_G}\dashv M^\psi$.
\end{enumerate}
\end{thm}

\begin{rem}
We can also define a category 
$\CC_G$
in which the objects are 
triples $(A,\AA,i)$ where 
$(\AA, i) $ is a $G$-grading of a \csta\ $A$ and the morphisms are defined as in Proposition \ref{catone}. It is not difficult to see that all our results in the present subsection concerning $\CC^{\rm top}_G$ can be reformulated with $\CC_G$ instead. We have chosen to concentrate our attention on topologically $G$-graded $C^*$-algebras because this seems necessary in order to be able to define the normalization functor in the next subsection.

\end{rem}

\section{The normalization functor}\label{normalization}

In this section we will apply our monad technique,
and for this we will need:
\begin{enumerate}
\item
a functor $N:\CC^{\rm top}_G\to \CC^{\rm top}_G$,
and

\item
a natural transformation $\eta:\id_{\CC^{\rm top}_G}\to N$
\end{enumerate}
such that
$N\psi$ and $\psi N$ are both natural isomorphisms from $N$ to $N^2$.

\begin{prop}
For every 
morphism $\varphi:\AA\to \BB$
there is a unique \hm\ $\varphi^n$ making the diagram
\[
\begin{tikzcd}
C^*(\AA) \arrow[r,"\varphi^m"] \arrow[d,"\Lambda_\AA"']
&C^*(\BB) \arrow[d,"\Lambda_\BB"]
\\
C^*_r(\AA) \arrow[r,"\varphi^n"',dashed]
&C^*_r(\BB)
\end{tikzcd}
\]
commute.
\end{prop}

\begin{proof}
The argument is quite similar to that of \propref{eta prop}.
Enlarge the diagram:
\begin{equation}\label{enlarge}
\begin{tikzcd}
C^*(\AA) \arrow[r,"\varphi^m"] \arrow[d,"\Lambda_\AA"'] \arrow[dr,"\rho"']
&C^*(\BB) \arrow[d,"\Lambda_\BB"]
\\
C^*_r(\AA) \arrow[r,"\varphi^n"',dashed] \arrow[d,"P_\AA"']
&C^*_r(\BB) \arrow[d,"P_\BB"]
\\
\Gamma_0(\AA^0) \arrow[r,"\varphi^0"']
&\Gamma_0(\BB^0),
\end{tikzcd}
\end{equation}
where we define $\rho$ to make the top triangle commute.
The outer rectangle of the diagram~\eqref{enlarge} commutes when at the top we restrict to
$i_\AA(\Gamma_c(\AA))\to i_\BB(\Gamma_c(\BB))$,
and hence commutes at the \cst-level by density and continuity.

It suffices to show that
\[
\ker \Lambda_\AA\subset \ker \rho.
\]
The important point is that the conditional expectation $P_\BB$ is faithful.
Thus, if $a\in \ker\Lambda_\AA$ then so is $a^*a$, and then following the path down the left-hand side of diagram~\eqref{enlarge} and along the bottom gives the same as $P_\BB\circ \rho$,
and since the conditional expectation $P_\BB$ is faithful we conclude that $\rho(a^*a)=0$, and so $a\in \ker\rho$.
\end{proof}

\begin{prop}\label{nor}
There is a functor $N:\CC^{\rm top}_G\to \CC^{\rm top}_G$ given on objects by
\[
N(A,\AA,i)=(\cstra,\AA,i_{\AA,r})
\]
and on morphisms $(\phi,\varphi):(A,\AA,i)\to (B,\BB,j)$ by
\[
N(\phi,\varphi)=(\varphi^n,\varphi):(\cstra,\AA,i_{\AA,r})\to (\cstrb,\BB,i_{\BB,r}).
\]
\end{prop}
\begin{proof} Straightforward.
\end{proof}

We now promote our \hm s $\eta_{(A, \AA, i)}:A\to \cstra$ 
to a natural transformation $\eta:\id_{\CC^{\rm top}_G}\to N$.
We must actually have a morphism in $\CC^{\rm top}_G$, so we embellish also here:
\[
(\eta_{(A,\AA,i)},\id_\AA):(A,\AA,i)\to N(A, \AA, i)
\]
and identify $\eta_{(A, \AA, i)}$ with $(\eta_{(A,\AA,i)},\id_\AA)$.

\begin{prop}\label{eta natural}
The \hm s $\eta_{(A,\AA,i)}$ combine to give a natural transformation
\[
\eta:{\rm id}_{\CC^{\rm top}_G}\to N.
\]
\end{prop}

\begin{proof}
This follows immediately from commutativity of the diagram
\[
\begin{tikzcd}
(A,\AA,i) \arrow[r,"{(\phi,\varphi)}"] \arrow[d,"\eta_{(A,\AA,i)}"']
&(B,\BB,j) \arrow[d,"\eta_{(B,\BB,j)}"]
\\
(\cstra,\AA,i_{\AA,r}) \arrow[r,"{(\varphi^n,\varphi)}"']
&(\cstrb,\BB,i_{\BB,r}).
\end{tikzcd}
\]
\end{proof}

Now we take care of the compositions $N\eta$ and $\eta M$.

\begin{prop}\label{N eta}
$N\eta:N\to N^2$ and $\eta N:N\to N^2$ are natural isomorphisms.
\end{prop}

\begin{proof}
We begin by computing $N^2$ and $N\eta$.

First, we claim that in fact $N^2=N$.
For an object $(A,\AA,i)$, we have
\[
N^2(A,\AA,i)
=N(\cstra,\AA,i_{\AA,r})
=(\cstra,\AA,i_{\AA,r})
=N(A,\AA,i).
\]
For a morphism $(\phi,\varphi):(A,\AA,i)\to (B,\BB,j)$ we have
\[
N(\phi,\varphi)
=(\varphi^n,\varphi):
(\cstra,\AA,i_{\AA,r})\to (\cstrb,\BB,i_{\BB,r})
\]
and normalizing this does nothing, so
\[
N^2(\phi,\varphi)
=N(\varphi^n,\varphi)
=(\varphi^n,\varphi)
=N(\phi,\varphi).
\]
Therefore we do have $N^2=N$. 

We turn to $N\eta$. For an object $(A,\AA,i)$,
we have
\[\eta_{(A,\AA,i)}=(\eta_{(A,\AA,i)},\id_\AA): (A,\AA,i) \to (\cstra,\AA,i_{\AA,r}),\]
so
\begin{align*}
N\eta_{(A,\AA,i)}
&=((\id_\AA)^n,\id_\AA)
=(\id_{C_r^*(\AA)},\id_\AA)
\\&=\id_{(\cstra,\AA,i_{\AA, r})}
=\id_{N(A,\AA,i)}.
\end{align*}
Therefore we get
$N\eta=\id:N\to N$, which is certainly a natural isomorphism.

Next we check that $\eta N:N\to N^2$ is a natural isomorphism.
For an object $(A,\AA,i)$ we have
$N(A,\AA,i)=(\cstra,\AA,i_{\AA,r})$,
and in 
Proposition \ref{eta prop} the
 relevant diagram is
\[
\begin{tikzcd}
\gam A \arrow[r,"i_{\AA,r}"] \arrow[d,"i_{\AA,r}"']
&\cstra \arrow[dl,"\eta_{(\cstra,\AA,i_{\AA,r})}"]
\\
\cstra,
\end{tikzcd}
\]
so by uniqueness we must have $\eta_{(\cstra,\AA,i_{\AA,r})}=\id_{\cstra}$. Therefore we have
\[
(\eta N)_{(A,\AA,i)}=\eta_{N(A,\AA,i)}=(\id_{C_r^*(\AA)},\id_{\AA})=\id_{(C_r^*(\AA),\AA,i_{\AA, r})},
\]
which is an isomorphism. Thus  $\eta N$ is a natural isomorphism.

\end{proof}

\begin{defn}
We denote by $\NN_G$
the full subcategory of $\CC^{\rm top}_G$ whose objects $(A,\AA,i)$
satisfy the property that
$\eta_{(A, \AA, i)}:A\to \cstra$
is an isomorphism. In other words, $\NN_G= (\CC^{\rm top}_G)^\eta$.
\end{defn}

We now immediately conclude the following from \propref{refl}:

\begin{thm}
Let $N:\CC^{\rm top}_G\to \CC^{\rm top}_G$ and $\eta:\id_{\CC^{\rm top}_G}\to N$ be as above.
Then:
\begin{enumerate}
\item
$(N,\eta)$ is an idempotent monad on $\CC^{\rm top}_G$.

\item
$\inc_{\NN_G}:\CC^{\rm top}_G\to \NN_G$ is left-adjointable and
$N^\eta$ is a left adjoint.

\item
$\NN_G$ is a reflective subcategory of $\CC^{\rm top}_G$.

\item
$\eta$ is the unit of the adjunction $N^\eta\dashv \inc_{\NN_G}$.
\end{enumerate}
\end{thm}

\section{Maximal-normal equivalence}\label{max-nor}
We consider now the category $\CC^{\rm top}_G$ of topologically graded $C^*$-algebras over a given locally compact Hausdorff \' etale goupoid $G$, equipped with the functors $M, N$ and the natural transformations $\psi, \eta$ constructed in the previous subsections. It is then clear that the setup in subsection  \ref{MandN} is satisfied. Our goal is to show that the full coreflective  subcategory $\MM_G= (\CC^{\rm top}_G)^\psi$ 
is  equivalent to the full reflective subcategory $\NN_G= (\CC^{\rm top}_G)^\eta$,
and that these subcategories provide a new example of a pair satisfying the ``maximal-normal" type equivalence introduced in \cite{reflective}, hence of an associated pair (\cite{KL, asspairs}).

\begin{prop} \label{MetaNpsi} The natural transformations  $M\eta: M\to MN$ and $N\psi: NM \to N$ are both isomorphisms.
\end{prop}
\begin{proof} Let $(A, \AA, i)$ be an object in $\CC^{\rm top}_G$ and set 
$(\phi, \varphi):= M\eta_{(A, \AA,i)}$. 
Then
\[M(A,\AA, i)= (C^*(\AA), \AA, i_\AA) = M(C_r^*(\AA), \AA, i_{\AA, r}),\]
so $(\phi, \varphi):(C^*(\AA), \AA, i_{\AA})\to(C^*(\AA), \AA, i_{\AA})$, 
and the following diagram
\[
\xymatrix{
(C^*(\AA), \AA, i_\AA) \ar[r]^-{\psi_{(A, \AA, i)}} \ar
[d]_{(\phi, \varphi)}
&(A, \AA, i) \ar[d]^{\eta_{(A, \AA, i)}}
\\
(C^*(\AA), \AA, i_\AA) \ar[r]_-{\psi_{(C_r^*(\AA), \AA, i_{\AA, r})}}
&(C_r^*(\AA), \AA, i_{\AA, r})
}
\]
is commutative, i.e., we have
\begin{equation}\label{Meta}
 \psi_{(C_r^*(\AA), \AA, i_{\AA,r})} \circ (\phi, \varphi) = \eta_{(A, \AA,i)} \circ \psi_{(A, \AA,i)}.
 \end{equation}
In fact, the universal property of Proposition \ref{aup} tells us that $(\phi, \varphi)$ is the \emph{unique} morphism in $\CC^{\rm top}_G$  such that (\ref{Meta}) holds, i.e., satisfying
\[ (\psi_{C_r^*(\AA)} \circ \phi, \varphi) = (\eta_A\circ \psi_A, {\rm id}_\AA), \]
that is, $\psi_{C_r^*(\AA)} \circ \phi = \eta_A\circ \psi_A$ and  $\varphi = {\rm id}_\AA$.

Now, $\psi_{C_r^*(\AA)}:C^*(\AA)\to C_r^*(\AA)$ is the unique homomorphism such that 
\[\psi_{C_r^*(\AA)}\circ i_{\AA} = i_{\AA, r}.\] This means that $\psi_{C_r^*(\AA)} = \Lambda_\AA= \eta_A\circ \psi_A$, giving that 
\[ 
 \psi_{(C_r^*(\AA), \AA, i_{\AA,r})} \circ ({\rm id}_{C^*(\AA)},  {\rm id}_\AA)  = 
  (\eta_A\circ \psi_A, {\rm id}_\AA).\]
Using the uniqueness property of $(\phi, \varphi)$, we get that $(\phi, \varphi) = ({\rm id}_{C^*(\AA)},  {\rm id}_\AA)$, hence that  $M\eta_{(A, \AA,i)}=(\phi, \varphi)$  is an isomorphism. This shows that the natural transformation $M\eta$ is an isomorphism. 

By a dual argumentation, we get that $N\psi_{(A, \AA,i)} = ({\rm id}_{C_r^*(\AA)},  {\rm id}_\AA)$, hence that $N\psi$ is also an isomorphism. 

\end{proof}

As recalled in subsection \ref{MandN}, an immediate consequence of Proposition \ref{MetaNpsi} is the following.
\begin{thm}\label{max nor grading} The categories $\MM_G$ and $\NN_G$ are equivalent. More precisely, the adjunction $N^\eta \circ {\rm Inc}_{\MM_G} \dashv M^\psi \circ {\rm Inc}_{\NN_G}$ is a ``maximal-normal" adjoint equivalence  between $\MM_G$ and $\NN_G$. 
We also have
\[N^\eta \simeq (N^\eta \circ {\rm Inc}_{\MM_G}) \circ M^\psi \, \text{ and }\, M^\psi \simeq (M^\psi \circ {\rm Inc}_{\NN_G}) \circ N^\eta.\]
\end{thm}  
 
%\vspace{-6ex}

\end{document}